\documentclass[12pt]{article}

\usepackage{amssymb,amsmath,mathtools,amsthm}
\usepackage{graphicx,float}
\usepackage{enumerate}
\usepackage{cite}
\usepackage{subfigure} 
\usepackage{xcolor}
\usepackage{geometry}

\usepackage[plainpages=false,hypertexnames=false,pdfpagelabels=true,
hyperindex=true,linktocpage,pagebackref=true,pdfa=true,hidelinks,colorlinks=true,citecolor=violet
,linkcolor=blue]{hyperref}

\newtheorem{theorem}{Theorem}[section]
\newtheorem{corollary}[theorem]{Corollary}
\newtheorem{lemma}[theorem]{Lemma}

\newtheorem{remark}{Remark}

\begin{document}

\title{Signless Laplacian determinations of some graphs with independent edges}

\author{R. Sharafdini$^b$\footnote{sharafdini@pgu.ac.ir}\qquad\qquad A.Z. Abdian$^a\footnote{abdian.al@fs.lu.ac.ir; aabdian67@gmail.com; azeydiabdi@gmail.com}$
\\[3mm]
$^b$ Department of Mathematics, Faculty of Science,\\
Persian Gulf University, Bushehr, 7516913817, Iran\\[5mm]
$^a$ Department of the mathematical Science, College of Science,\\
 Lorestan University, Lorestan, Khoramabad 41566, Iran
}

\maketitle

\begin{abstract}
  {Signless Laplacian determinations of some graphs with independent edges}%
{Let $G$ be a simple undirected graph. Then the signless Laplacian matrix of $G$ is
defined as $D_G + A_G$ in which $D_G$ and $A_G$ denote the degree matrix and the adjacency matrix of $G$,
respectively. The graph $G$ is said to be determined by its signless Laplacian spectrum ({\rm DQS}, for short), if any graph
having the same signless Laplacian spectrum as $G$ is isomorphic to $G$.
We show that $G\sqcup rK_2$ is determined by its signless Laplacian spectra
under certain conditions, where $r$ and $K_2$ denote a natural number and the complete graph on two vertices, respectively. Applying
these results, some   {\rm DQS}   graphs with independent edges are obtained.}\\[1cm]
\textbf{2010 Mathematics Subject Classification:}05C50\\
\textbf{Keywords}:Spectral characterization, Signless Laplacian spectrum, Cospectral graphs, Union of graphs
\end{abstract}

\section{Introduction}
All graphs considered here are simple and undirected. All notions on graphs that are not defined here can be found in \cite{CRS10,BroHamSpec}.
Let $G$ be a simple graph with the vertex set $ V = V (G) = \left\{ {v_1,\ldots , v_n} \right\}$ and the edge set $E = E(G)$. Denote by $d_i$ the degree of the vertex $ v_i $. The \emph{adjacency matrix}
$A_G$ of $G$ is a square matrix of order $n$, whose $(i,j)$-entry
is 1 if $v_i$ and $v_j$ are adjacent in $G$ and $0$
otherwise. The \emph{degree matrix} $D_G$ of $G$ is a diagonal matrix of order $n$ defined as $D_G=\mathrm{diag}(d_1,\ldots,d_n)$. The matrices $ L_G=D_G-A_G $ and $ Q_G=D_G+A_G $ are called the \emph{Laplacian matrix} and the \emph{signless Laplacian matrix} of $G$, respectively.
The multiset of eigenvalues of $Q_G$ (resp.
$L_G$, $A_G$) is called the \emph{$Q$-spectrum} (resp. \emph{$L$-spectrum}, \emph{$A$-spectrum}) of $G$. For any bipartite graph, its $Q$-spectrum
coincides with its $L$-spectrum. Two graphs are $Q$-cospectral (resp. $L$-cospectral, $A$-cospectral) if they have the
same $Q$-spectrum (resp. $L$-spectrum, $A$-spectrum). A graph $G$ is said to be   {\rm DQS}   (resp.  {\rm DLS},  {\rm DAS}) if there is
no other non-isomorphic graph $Q$-cospectral (resp. $L$-cospectral, $A$-cospectral) with $G$. Let us denote the $Q$-spectrum of $G$ by ${\rm Spec}_{Q}(G)= \left\{ {[q_1]^{ m_1 },[q_2]^{ m_2 },\ldots,[q_n]^{m_n }}\right\}$, where $ m_i $ denotes the multiplicity of $ q_i $ and $q_1\ge q_2\ge \ldots \ge q_n$.

The \emph{join} of two graphs  $G$  and $ H $ is a graph formed from disjoint copies of  $G$  and $ H $ by connecting each vertex of  $G$  to each vertex of $H$. We denote the join of two graphs  $G$  and $ H $ by $ G\nabla H $. The complement of a graph $G$  is denoted by $ \overline{G} $. For two disjoint graphs $G$ and $H$, let $G\sqcup H$ denotes
the disjoint union of $G$ and $H$, and $rG$ denotes the disjoint union of $r$ copies of $G$, i.e.,  $rG=\underbrace{G\sqcup\ldots \sqcup G}_{r-\text{times}}$.

Let $G$ be a  connected graph  with $n$ vertices and $m$ edges. Then $G$ is called \emph{unicyclic} (resp. \emph{bicyclic}) if $m=n$ (resp.
$m = n + 1$). If $G$ is a unicyclic graph containing an odd (resp. even) cycle, then $G$ is called \emph{odd unicyclic} (resp. \emph{even unicyclic}).

Let $C_n$, $P_n$, $K_n$ be the cycle, the path and the complete graph of order $n$, respectively.
$K_{s,t}$ the complete bipartite graph with $s$ vertices in one part
and $t$ in the other.

Let us remind that the \emph{coalescence}\cite{CvetDoobSach1995} of two graphs $G_1$ with distinguished vertex $v_1$ and $G_2$ with distinguished vertex $v_2$, is formed by identifying vertices $v_1$ and $v_2$ that is, the vertices $v_1$ and $v_2$ are replaced by a single vertex $v$ adjacent to the same vertices in $G_1$ as $v_1$ and the same vertices in $G_2$ as $v_2$. If it is not necessary $v_1$ or $v_2$ may not be specified.

The \emph{friendship graph} $F_n$ is a graph with $2n+1$ vertices and $3n$ edges obtained by
the coalescence of $n$ copies of $C_3$ with a common vertex as the distinguished vertex; in fact, $F_n$ is nothing but $K_1\nabla nK_2$.

\begin{figure}[H]
\centering
\subfigure[{}]
{\label{fig:T}\includegraphics[width=4cm]{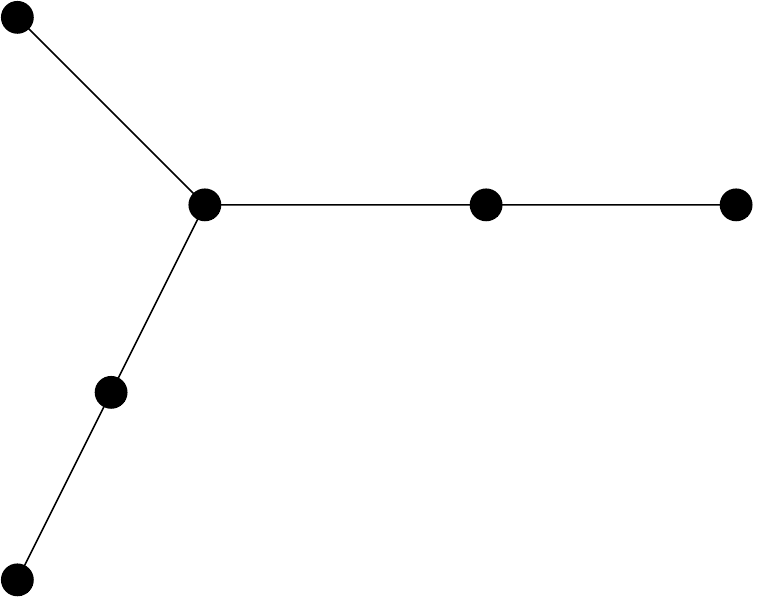}}
\quad\quad\quad\quad\quad\quad\quad\quad
\subfigure[{}]
{\label{fig:lollipop}\includegraphics[width=5cm] {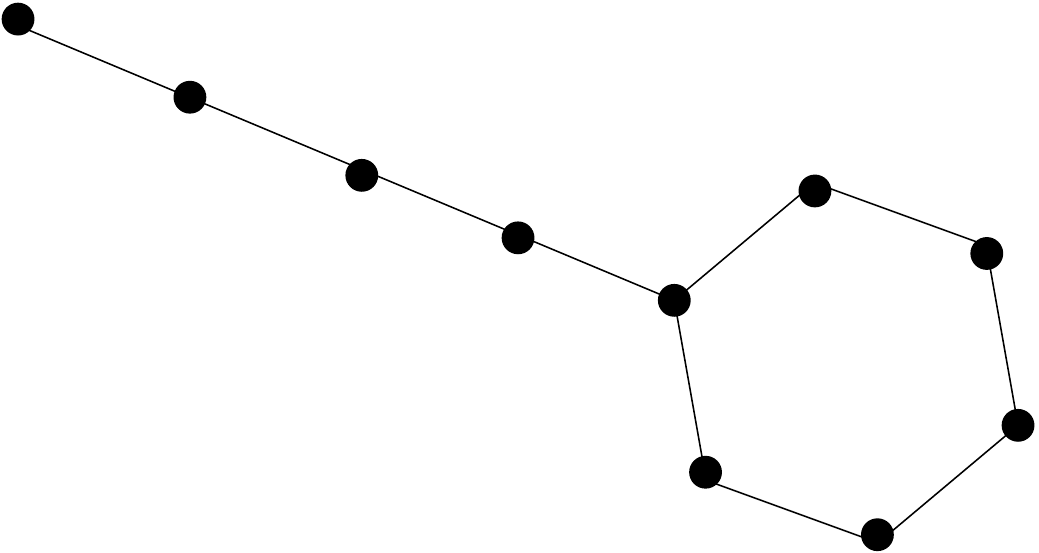}}
\caption{(a)The T-shape tree $T(0,1,1)$,\quad(b)The lollipop graph $H_{11,6}$} \label{fig:T-lollipop}
\end{figure}

The \emph{lollipop graph},
denoted by $H_{n,p}$, is the
coalescence of a cycle $C_p$ with arbitrary distinguished vertex and a path $P_{n-p}$ with a pendent vertex
as the distinguished vertex; for example $H_{11,6}$ is depicted in Figure \ref{fig:lollipop}.
We denote by $T(a, b, c)$ the \emph{$T$-shape} tree obtained by identifying the end vertices of three paths $P_{a+2}$, $P_{b+2}$ and $P_{c+2}$. In fact, $T(a, b, c)$ is a tree with one and only one vertex $v$ of degree $3$ such that $T(a, b, c)-\{v\}=P_{a+1}\sqcup P_{b+1}\sqcup P_{c+1}$; for example $T(0,1,1)$ is depicted in Figure \ref{fig:T}.

\begin{figure}[H]
  \centering
  \includegraphics[width=5cm]{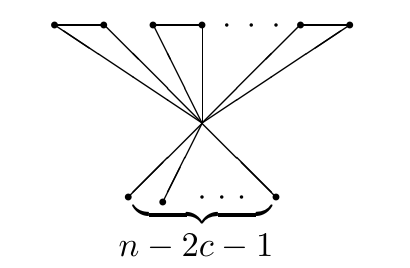}
  \caption{The bundle $S(n,c)$.}\label{fig:snc}
\end{figure}

A \emph{kite graph} $Ki_{n,w}$ is a graph obtained from a clique $K_w$ and a path
$P_{n-w}$
is the coalescence of $K_w$ with an arbitrary distinguished vertex and a path $P_{n-w+1}$ with a pendent vertex
as the distinguished vertex. Note that $K_{5,3}$ and $K_{6,3}$ are identical.
A tree is called \emph{starlike} if it has exactly one vertex of degree greater than two. We denote by $U_{r,n-r}$ the graph obtained by attaching $n-r$ pendent vertices to a vertex of $C_r$. In fact, $U_{r,n-r}$
is the coalescence of $K_{1,n-r-1}$ and $P_{n-w+1}$ where distinguished vertices are
the vertex of degree $n-r$ and a pendent vertex, respectively.
A graph is a \emph{cactus}, or a \emph{treelike graph}, if any pair of its cycles has at most one common vertex \cite{Radosavl}. If all cycles of the cactus $G$ have exactly one common vertex, then $G$ is called a \emph{bundle} \cite{Borov}. Let $S(n,c)$ be the bundle with $n$ vertices and $c$ cycles of length 3 depicted in Figure \ref {fig:snc}, where
$n\ge2c+1$ and $c\ge0$. By the definition, it follows that
$S(n,c) = K_1\nabla(cK_2\sqcup (n-2c-1)K_1)$. In fact $S(n,c)$
is the coalescence of $F_c$ and $K_{1,n-2c-1}$ where the distinguished vertices are
the vertex of the degree $2c$ and the vertex of the degree $n-2c-1$, respectively.

\begin{figure}[H]
  \centering
  \includegraphics[width=4cm]{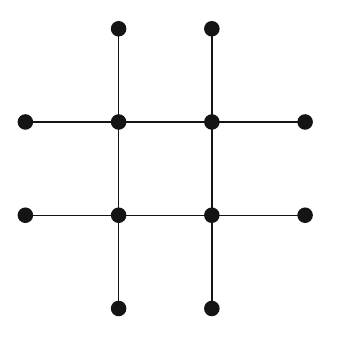}
  \caption{$C_4\circ 2K_1$.}\label{fig:C42K1}
\end{figure}

Let $G$ be a graph with $n$ vertices, $H$ be a graph with $m$ vertices. The \emph{corona} of
$G$ and $H$, denoted by $G\circ H$, is the graph with $n+mn$ vertices obtained from $G$
and $n$ copies of $H$ by joining the $i$-th vertex of $G$ to each vertex in the $i$-th copy of $H~(i\in\{1,\ldots, n\})$; for example $C_4\circ 2K_1$ is depicted in Figure \ref{fig:C42K1}.

\begin{figure}[H]
  \centering
  \includegraphics[width=9cm]{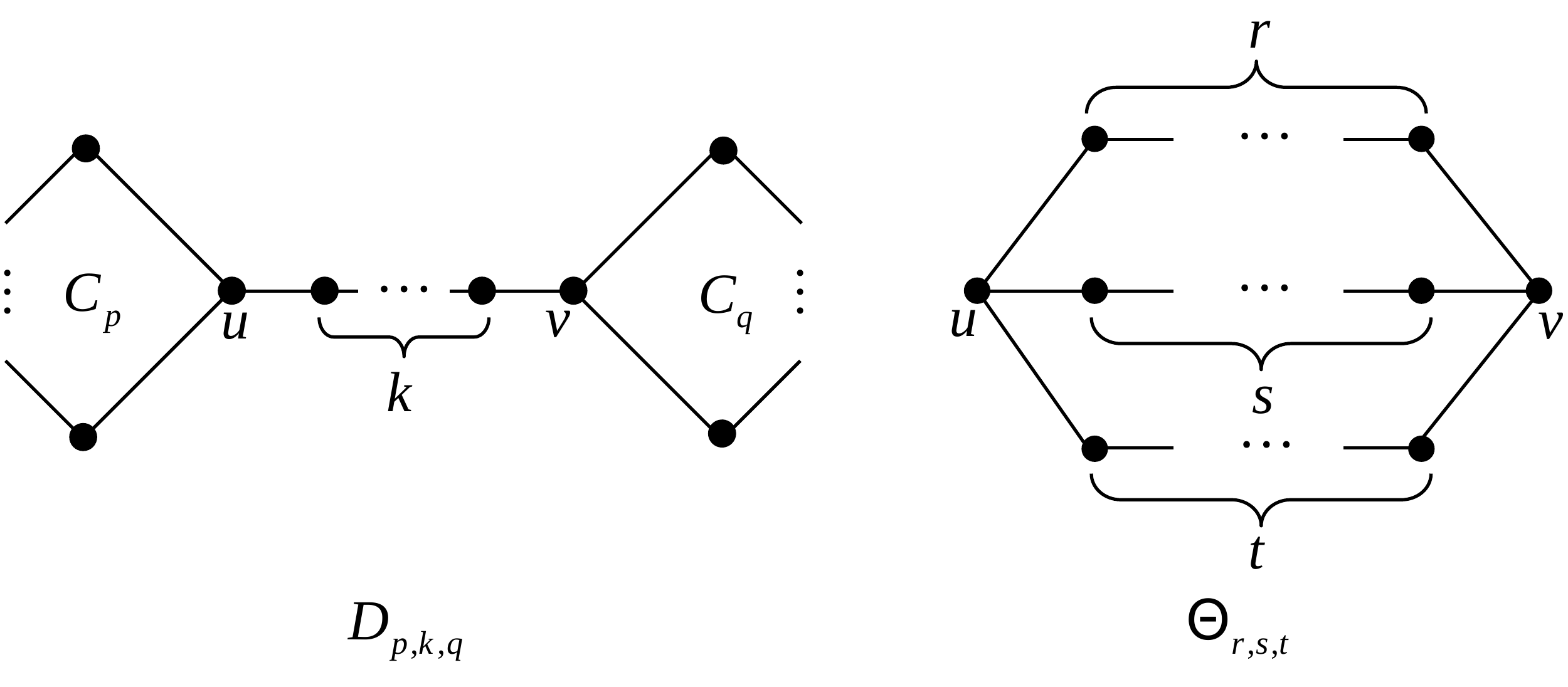}
  \caption{The graphs $D_{p,k,q}$ and $\Theta_{r,s,t}$.}\label{fig:DumTheta} 
\end{figure}

A \emph{complete split graph} $CS(n,\alpha)$, is a graph on $n$ vertices consisting of a clique on $n-\alpha$ vertices and an independent set on the remaining $\alpha~(1\leq\alpha\leq n-1)$ vertices in which each vertex of the clique is adjacent to each vertex of the independent set.
The \emph{dumbbell graph}, denoted by $D_{p,k,q}$, is a bicyclic graph obtained from two cycles $C_p$, $C_q$ and a path $P_{k+2}$ by identifying each pendant vertex of $P_{k+2}$ with a vertex of a cycle, respectively. The \emph{theta graph}, denoted by $\Theta_{r,s,t}$, is the graph formed by joining two given vertices via three disjoint paths $P_{r}$, $P_{s}$ and $P_{t}$, respectively, see Figure \ref{fig:DumTheta}.

The problem "which graphs are determined by their spectrum?" was posed by G\"{u}nthard and Primas \cite{Gunthard} more than
60 years ago in the context of H\"{u}ckel's theory in chemistry.
In the most recent years mathematicians have  devoted their attention to this problem and many papers focusing on this topic are now appearing.
In \cite{VH} van Dam and Haemers  conjectured that almost all graphs are determined by their spectra.
Nevertheless, the set of graphs that are known to be determined by their spectra is too small. So, discovering infinite classes of graphs that are determined by their spectra can be an interesting problem.
Cvetkovi\'{c}, Rowlinson and Simi\'{c} in \cite{CRS1,CS,CSI,CSIII} discussed the development of a spectral theory of graphs
based on the signless Laplacian matrix, and gave several reasons why it is superior to other graph matrices such as the adjacency and the Laplacian matrix. It is interesting to construct new {\rm DQS}   ({\rm DLS}) graphs from known   {\rm DQS}   ({\rm DLS}) graphs. Up to now, only some graphs with special structures are shown to be  {\it determined by their spectra} (DS, for short)
(see \cite{A, AA,AAA, AAAA, AAAAA,AAAAAA,AAAAA1,CRS1,CS,LiuLuEJLA,M,Mer
,MK,CHA,LC,HuangK,Liu2,ZhagLollipo,Das-Kite,Wang-FR,WangCMJ,Omidi,Liu3,Das-CSG,LiuLiuWei,WangDumbbell}
and the references cited in them).
About the background of the question
"Which graphs are determined by their spectrum?", we refer to \cite{VH},\cite{vandam2009}.
For a  {\rm DQS} graph $G$,  $G\nabla K_2$ is also  {\rm DQS} under some conditions \cite{LC}. A graph is  {\rm DLS} if
and only if its complement is  {\rm DLS}. Hence we can obtain  {\rm DLS} graphs from known  {\rm DLS} graphs by adding independent edges.
In \cite{HuangK} it was shown that $G\sqcup rK_1$ is   {\rm DQS}   under certain conditions.
In this paper, we investigate signless Laplacian spectral characterization of graphs with independent edges. For a   {\rm DQS}   graph $G$, we show that $G \sqcup rK_2$ is {\rm DQS}   under certain conditions. Applying these results, some   {\rm DQS}   graphs with independent edges are obtained.\\

\section{Preliminaries}
In this section, we give some lemmas which are used to prove our main results.

\begin{lemma}[\cite{CS,CRS1}]\label{lem 2-1}
Let $G$ be a graph.
For the adjacency matrix of $G$, the
following can be deduced from the spectrum:
\begin{enumerate}
\item[{\rm(1)}]  The number of vertices.
\item[{\rm(2)}]  The number of edges.
\item[{\rm(3)}]  Whether $G$ is regular.
\end{enumerate}
For the Laplacian matrix, the following follows from the spectrum:
\begin{enumerate}
\item[{\rm(4)}] The number of components.
\end{enumerate}
For the signless Laplacian matrix, the following follow from the spectrum:
\begin{enumerate}
\item[{\rm(5)}] The number of bipartite components, i.e., the multiplicity of the eigenvalue 0 of the signless Laplacian matrix is equal to the
number of bipartite components.
\item[{\rm(6)}] The sum of the squares of degrees of vertices.
\end{enumerate}
\end{lemma}
\begin{lemma}[\cite{CRS1}]\label{lem:gipower}
Let  $G$  be a graph with $n$ vertices, $ m $ edges, $ t $ triangles and the vertex degrees $ d_1,d_2,\ldots,d_n$. If $ {T_k} = \sum\limits_{i = 1}^n q_i(G)^k$, then we have
\[T_0=n,\quad T_1=\sum\limits_{i = 1}^n d_i=2m,\quad T_2=2m+\sum\limits_{i = 1}^n d^2_i,\quad  T_3=6t+3\sum\limits_{i = 1}^n d^2_i+\sum\limits_{i = 1}^n d^3_i.\]
\end{lemma}

For a graph $G$, let $P_L(G)$ and $P_Q(G)$ denote the product of all nonzero eigenvalues of $L_G$ and $Q_G$, respectively. Note that $P_L(K_2) = P_Q(K_2) = 2$.
We assume that $P_L(G) = P_Q(G) = 1$ if $G$ has no edges.

\begin{lemma}[\cite{CRS10}]\label{lem 2-3}  For any connected bipartite graph $G$ of order $n$, we have $P_Q(G) = P_L(G) = n\tau(G)$, where $\tau(G)$ is the
number of spanning trees of $G$. Especially, if $T$ is a tree of order $n$, then $P_Q(T)= P_L(T)= n$.
\end{lemma}

\begin{lemma}[\cite{MK}]\label{lem 2-4} Let $G$ be a graph with $n$ vertices and $m$ edges.
\begin{enumerate}[\rm (i)]
  \item $\det(Q_G) = 4$ if and only if $G$ is an odd unicyclic graph.
  \item  If $G$ is a non-bipartite connected graph and $m>n$, then $\det(Q_G)\ge 16$, with equality if and only if  $G$  is a non-bipartite bicyclic graph with $ C_4 $ as its induced subgraph.
\end{enumerate}
\end{lemma}

\begin{lemma}[\cite{CRS10}]\label{lem:q1GH}
Let $e$ be any edge of a graph $G$ of order $n$. Then
\[q_1(G) \ge q_1(G-e) \ge q_2(G) \ge q_2(G-e) \ge\ldots\ge q_n(G) \ge q_n(G-e) \ge 0.\]
\end{lemma}

\begin{lemma}[\cite{CvetDoobSach1995}]\label{lem:subg-q1} Let $H$ be a proper subgraph of a connected graph $G$. Then $q_1(G) > q_1(H)$.
\end{lemma}

\begin{lemma}[\cite{CvetDoobSach1995}]\label{lem:q1}
Let $G$ be a graph with $n$ vertices and $m$ edges. Then
$q_1(G) \ge \frac{4m}{n}$, with equality if and only if $G$ is
regular.
\end{lemma}

\begin{lemma}[\cite{CRS1}]\label{lem:Pathq1}
For a graph $G$,
$0<q_1(G)< 4$ if and only if all components of $G$ are paths.
\end{lemma}

\begin{lemma}[\cite{VH}]\label{lem:regularDAQS}
A regular graph is {\rm DQS} if and only if it is  {\rm DAS}.
A regular graph $G$ is {\rm DAS} {\rm (DQS)} if and only if $\overline{G}$ is {\rm DAS} {\rm (DQS)}.
\end{lemma}

\begin{lemma}[\cite{CS}]\label{lem:regularDAS}
Let $G$ be a $k$-regular graph of order $n$. Then $G$ is  {\rm DAS}  when $k\in\{0, 1, 2, n-3, n-2, n-1\}$.
\end{lemma}

\begin{lemma}[\cite{CHA}]\label{lem:RegularCone}
Let $G$ be a $k$-regular graph of order $n$. Then $G\nabla K_1$ is
 {\rm DQS} for $k\in\{1, n-2\}$, for $k = 2$ and $n\ge 11$. For $k= n-3$,
$G\nabla K_1$ is  {\rm DQS} if and only if $\overline{G}$ has no triangles.
\end{lemma}

\begin{lemma}[\cite{LC}]\label{lem:RegularConeK2}
Let $G$ be a $k$-regular graph of order $n$. Then $G\nabla K_2$ is
 {\rm DQS} for $k\in\{1,n-2\}$. For $k= n-3$,
$G\nabla K_2$ is  {\rm DQS} if and only if $\overline{G}$ has no triangles.
\end{lemma}

\begin{lemma}[\cite{HuangK}]\label{lem:KnMatching}
The following hold for graphs with isolated vertices:
\begin{enumerate}[{\rm(i)}]
\item
Let $T$ be a {\rm DLS} tree of order $n$. Then $T\sqcup rK_1$ is {\rm DLS}. If $n$ is not divisible by 4, then $T\sqcup rK_1$ is {\rm DQS}.
\item
The graphs $\overline{P_n}$ and $\overline{P_n}\sqcup rK_1$ are {\rm DQS}.
\item
Let $G$ be a graph obtained from $K_n$ by deleting a matching. Then
$G$ and $G\sqcup rK_1$
 are {\rm DQS}.
\item
A ($n-4$)-regular graph of order $n$ is
{\rm DAS} {\rm (DQS)} if and only if its complement is a 3-regular {\rm DAS} {\rm (DQS)} graph.
\item
Let $G$ be a ($n-3$)-regular graph of order $n$. Then $G\sqcup rK_1$ is {\rm DQS}.
\end{enumerate}
\end{lemma}

Now let us list some known families of {\rm DQS} graphs.
\begin{lemma}\label{lem:main}
The following graphs are   {\rm DQS}.
\begin{enumerate}[{\rm(i)}]
\item
 The graphs $P_n$, $C_n$, $K_n$, $K_{m,m}$,
 $rK_{n}$,
 $P_{n_1} \sqcup P_{n_2}\sqcup\ldots \sqcup P_{n_k}$ and
 $C_{n_1} \sqcup C_{n_2}\sqcup\ldots \sqcup C_{n_k}$, \cite{VH}.
\item
Any wheel graph $K_1\nabla C_n$, \cite{Liu2}.
\item Every lollipop graph $H_{n,p}$, \cite{ZhagLollipo}.
\item Every kite graph $Ki_{n,n-1}$ for $n\ge 4$ and $n\neq 5$, \cite{Das-Kite}.
\item The friendship graph $ F_n $, \cite{Wang-FR}.
\item
$(C_n\circ tK_1)$, for $(n \neq 32, 64)$  and $t\in\{1,2\}$, \cite{BuCorona-GC},\cite{MK}.
\item
The line graph of a $T$-shape tree $T(a, b, c)$ except $T(t, t, 2t+1)~(t>1)$, \cite{WangCMJ}.
\item
The starlike tree with maximum degree 4,\cite{Omidi}.
\item $U_{r,n-r}$  for $r\ge 3$, \cite{Liu3}.
\item
$CS(n,\alpha)$ when
$1\leq\alpha\leq n-1$ and $\alpha\neq 3$, \cite{Das-CSG}.
\item
For $n\ge 2c+1$ and $c\ge 0$, $\overline{S(n,c)}$ and  $S(n,c)$ except for the case of $c = 0$ and $n=4$, \cite{LiuLiuWei}.
\item
$K_{1,n-1}$ for $n\neq 4$, \cite{LiuLiuWei}.
\item
$G\nabla K_m$ where $G$ is an
$(n-2)$-regular graph on $n$ vertices, and $\overline{K_n}\nabla K_2$ except for $n=3$, \cite{LiuLuEJLA}.
\item
All dumbbell graphs different from $D_{3q,0,q}$ and all theta graphs, \cite{WangDumbbell}.
\end{enumerate}
\end{lemma}

It is easy to see that $K_{1,3}$ and $K_3\sqcup K_1$ are $Q$-cospectral,
i.e., ${\rm Spec}_{Q}(K_{1,3})={\rm Spec}_{Q}(K_{3})=\{[4]^{1},[1]^{2},[0]^{1}\}$. Therefore,
$S(n, c)$ is not {\rm DQS} when $c=0$ and $n = 4$, since $S(n,0)$ is nothing but $K_{1,n-1}$.

\section{Main Results}
We first investigate spectral characterizations of the union of a tree and several complete graphs $ K_2 $.
\begin{theorem}\label{the 3-1}
 Let $T$ be a  {\rm DLS} tree of order $n$. Then $T\sqcup rK_2$ is  {\rm DLS} for any positive integer $r$. Moreover, if $n$ is odd and $r=1$, then $T\sqcup rK_2$ is {\rm DQS}.
 \end{theorem}

\begin{proof} For $ n, r\in\left\{ {1,2} \right\}$ see Lemma \ref{lem:KnMatching} (i) and Lemma \ref{lem:main} (i). So, one may suppose that $ n,r\geq 3 $. Let $G$ be any graph $L$-cospectral with $T\sqcup rK_2$. By Lemma \ref {lem 2-1}, $G$ has $n + 2r$ vertices, $n -1+r$ edges and
$r + 1$ components. So each component of $G$ is a tree. Suppose that $G = G_0 \sqcup G_1 \sqcup\ldots \sqcup  G_r$, where $G_i$ is a
tree with $n_i$ vertices and $n_0 \geq n_1 \geq \ldots \geq n_r \geq 2$. For $ n_i, n_ r\in\left\{ {1} \right\}$ see Lemma \ref{lem:KnMatching} (i) and Lemma \ref{lem:main} (i). Hence we consider $ n, n_i, r\geq 2 $.  Since $G$ is $L$-cospectral with $T\sqcup rK_2$, by Lemma \ref {lem 2-3}, we
get \[n_0n_1\ldots n_r=P_L(G_0)\ldots P_L(G_r)=P_L(G_0\sqcup\ldots \sqcup G_r)=P_L(G)=P_L(T)P_L(K_2)^r = n2^r.\].

We claim that $ n_r=2 $. Suppose not and so $ n_r\geq 3 $. This means that $n_0\geq n_1\geq ...\geq n_r\geq 3 $. Hence $n2^r =n_0n_1...n_r\geq 3^{r+1}$ or $ n(\dfrac{2}{3})^r\geq3 $. Now, if $r \to \infty$, then $ 0\geq 3 $, a contradiction. So, we must have $ n_r=2 $. By a similar argument one can show that $ n_1=...=n_{r-1}=2 $ and so $ n_0=n $. Hence $G = G_0 \sqcup rK_2$. Since $G$ and $T\sqcup rK_2$ are $L$-cospectral, $G_0$ and $T$ are
$L$-cospectral. Since $T$ is  {\rm DLS}, we have $G_0 = T$, and thus $G = T\sqcup rK_2$. Hence $T\sqcup rK_2$ is  {\rm DLS}.

Let $H$ be any graph $Q$-cospectral with $T\sqcup rK_2$. By Lemma \ref{lem 2-1}, $H$ has $n +2r$ vertices, $n -1+r$ edges and $r + 1$
bipartite components. So one of the following holds:
\begin{enumerate}
  \item[\rm(i)]  $H$ has exactly $r + 1$ components, and each component of $H$ is a tree.
  \item[\rm(ii)] $H$ has $r + 1$ components which are trees, the other components of $H$ are odd unicyclic.
\end{enumerate}
In what follows we show that (ii) does not occur if $n$ is odd and $r=1$.
If (ii) holds, then by Lemma \ref {lem 2-4}, $P_Q(H)$ is divisible by 4 since $H$ has a cycle of odd order as a component. Since $T$ is a tree of order $n$, by
Lemma \ref {lem 2-3}, $P_Q(H)=P_Q(T)P_Q(K_2)^r= n2^r$ is divisible by 4, a contradiction.
%
%
%
Therefore (i) must hold. In this case, $H$ and $T\sqcup rK_2$ are both bipartite, and so they are also $L$-cospectral. By the previous part, $T\sqcup rK_2$ is  {\rm DLS}. So we
have $H = T\sqcup rK_2$.

 Hence $T\sqcup rK_2$ is   {\rm DQS}   when $n$ is odd and $r=1$.
 \end{proof}
\begin{remark} Some  {\rm DLS} trees are given in \cite{HuangK} and references therein. We can obtain  some {\rm DLS} ({\rm DQS}) trees with independent edges from Theorem \ref{the 3-1}.
\end{remark}

Lemma \ref{lem:main} and Theorem \ref{the 3-1} imply the following corollary.

\begin{corollary} For an odd positive integer $n$, we have the following
\begin{enumerate}[{\rm (i)}]
\item
Let $T$ be a starlike tree of order $n$ and with maximum degree 4. Then $T\sqcup K_2$ is {\rm DQS}.
\item $P_n\sqcup K_2$ is {\rm DQS}.
\item
For $n\neq 4$, $K_{1,n-1}\sqcup K_2$  is {\rm DQS}.
\item
Let $\mathcal{L}$ be the line graph of a $T$-shape tree $T(a, b, c)$ except $T(t, t, 2t+1)~(t>1)$. Then $\mathcal{L} \sqcup K_2$ is {\rm DQS} if $a+b+c-3$ is odd.
\end{enumerate}
\end{corollary}


\begin{theorem}\label{the 3-2}
 Let $G$ be a   {\rm DQS} odd unicyclic graph of order $n\geq 7$. Then $G \sqcup rK_2$ is  {\rm DQS} for any positive integer $r$.
 \end{theorem}

\begin{proof}  Let $H$
be any graph $Q$-cospectral with $G \sqcup rK_2$. By Lemma \ref{lem 2-1}(5), 0 is not an eigenvalue of $G$ since it is an odd unicyclic. So by Lemma \ref{lem 2-4}, we have
$4=\det(Q_G)=P_Q(G)$. Moreover,
$$P_Q(H)=P_Q(G \sqcup rK_2)=P_Q(G)P_Q(K_2)^r=\det(Q_G)2^r=4.2^r=2^{r+2}.$$
By Lemma \ref {lem 2-1}, $H$ has $n + 2r$
vertices, $n+r$ edges and $r$ bipartite components. So one of the following holds:
\begin{enumerate}
  \item[{\rm (i)}]$H$ has exactly $r$ components each of which is a tree.
  \item [{\rm (ii)}] $H$ has $r$ components which are trees, the other components of $H$ are odd unicyclic.
\end{enumerate}
We claim that (i) does not hold, otherwise, we may assume that $H = H_1 \sqcup\ldots \sqcup H_r$, where $H_i$ is a tree with $n_i$ vertices and $n_1 \ge \ldots \ge n_r \ge 1$. It follows from Lemma \ref {lem 2-3} that
\[n_1\ldots n_r=P_Q(H_1)\ldots P_Q(H_r)=P_Q(H)=4.2^r=2^{r+2}.\]
So $n_1\ldots n_r=2^{r+2}$, $n_1\leq 8 $.
Since $G$ contains a cycle, say $C$, by Lemma \ref{lem:q1} we have
\begin{equation}\label{eq:HGC4}
q_1(H)=q_1(G)\ge q_1(C)=4.
\end{equation}
Let $\Delta(H)$ be the maximum degree of $H$. If $\Delta(H)\leq 2$, then all components of $H$ are paths, hence by Lemma \ref{lem:Pathq1},
$q_1(H)<4$, contradicting Eq. \eqref{eq:HGC4}. So $\Delta(H)\ge 3$. From $n_1\leq 8$ and $n_1\ldots n_r = 4.2^r=2^{(r+2)}$, we may assume that $H_1 = K_{1,7}$,
$H_2 =\ldots = H_r = K_2$. Since $H = K_{1,7} \sqcup (r-1)K_2$ has $n + 2r$ vertices, we get $n = 6$, a contradiction to $n\ge 7$.\\
If (ii) holds, then we may assume that $H = U_1\sqcup\ldots \sqcup U_c \sqcup H_1\sqcup\ldots \sqcup H_r$, where $U_i$ is odd unicyclic, $H_i$ is a tree with
$n_i$ vertices. By Lemmas \ref {lem 2-3} and \ref{lem 2-4}, $4.2^r= P_Q(H) = 4^cn_1\ldots  n_r$. So $c = 1$, $H_1 =\ldots = H_r = K_2$. Since
$H = U_1\sqcup rK_2$ and $G \sqcup rK_2$ are $Q$-cospectral, $U_1$ and $G$ are $Q$-cospectral. Since $G$ is   {\rm DQS}, we have $U_1 = G$,
$H = G \sqcup rK_2$.
 \end{proof}

\begin{remark}
Note that $C_4\sqcup 2P_3$ and $C_6\sqcup 2K_2$ are $Q$-cospectral, i.e.,
${\rm Spec}_{Q}(C_4\cup 2P_3)={\rm Spec}_{Q}(C_6\cup 2K_2)=
\{[4],[3]^2,[2]^2,[1]^2,[0]^3\}$. It follows that the condition "odd unicyclic of order $n\ge 7$" is essential in Theorem \ref{the 3-3}.
\end{remark}
\begin{remark} Some   {\rm DQS}   unicyclic graphs are given in
\cite{HuangK} and references therein. We can obtain  some {\rm DQS}   graphs with independent edges from Theorem \ref{the 3-2}.
\end{remark}
\begin{theorem}\label{the 3-3}
Let $G$ be a  {\rm DQS} graph of order $n\geq 5 $. If $G$  is non-bipartite  bicyclic graph with $C_4$ as its induced subgraph, then $G\sqcup rK_2$ is {\rm DQS} for any positive integer $r$.
\end{theorem}

\begin{proof} Let $H$ be any graph $Q$-cospectral with $G\sqcup rK_2$. By Lemma \ref {lem 2-4}, we have
\[P_Q(H) =P_Q(G\sqcup rK_2)=P_Q(G)P_Q(K_2)^r=P_Q(G)2^r\]

By Lemma \ref{lem 2-1}(5), 0 is not an eigenvalue of $G$ since it is non-bipartite.
So by Lemma \ref{lem 2-4}, we have
$16=\det(G_Q)=P_Q(G)$ and thus $P_Q(H)=16.2^r$.

By
Lemma \ref {lem 2-1}, $H$ has $n + 2r$ vertices, $n + 1+r$ edges and $r$ bipartite components. So $H$ has at least
$r - 1$ components which are trees.
Suppose that $H_1,H_2,\ldots , H_r$ are $r$ bipartite components of $H$, where $H_2,\ldots , H_r$ are trees. If $H_1$ contains an
even cycle, then by Lemmas \ref {lem 2-4} and \ref{lem:q1GH}, we have $P_Q(H) \geq P_Q(H_1) \geq 16$, and $P_Q(H) = 16.(2^{r-1})=2^{r-3}$ if and only if $H = C_4\sqcup (r-1)K_2$.
By $P_Q(H) = 16.(2^{r-1})=2^{r-3}$, we have $H = C_4\sqcup (r-1)K_2$. Since $H$ has $n + 2r$ vertices, we get $n = 2$, a contradiction ($G$
contains $C_4$). Hence $H_1,H_2,\ldots , H_r$ are trees.
Since $H$ has $n +2r$ vertices, $n + 1+r$ edges and $r$ bipartite components, $H$ has a non-bipartite component
$H_0$ which is a bicyclic graph. Lemmas \ref{lem 2-4} and \ref{lem:q1GH} imply that $P_Q(H)\geq P_Q(H_0)\geq 16$, and $P_Q(H) = 16.2^r$ if and only if
$H = H_0 \sqcup rK_2$ and $H_0$ contains $C_4$ as its induced subgraph. By $P_Q(H) = 16.2^r$, we have $H = H_0 \sqcup rK_2$. Since $H$
and $G \sqcup rK_2$ are $Q$-cospectral, $H_0$ and $G$ are $Q$-cospectral. Taking into account that $G$ is {\rm DQS}, we conclude that $H_0 = G, H = G \sqcup rK_2$.
Hence $G \sqcup rK_2$ is   {\rm DQS}.
\end{proof}

\begin{remark}Some   {\rm DQS}   bicyclic graphs are given in \cite{HuangK} and references therein. We can obtain   {\rm DQS}   graphs with independent
edges from Theorem \ref{the 3-3}.
\end{remark}

\begin{lemma}\label{lem:thm4}
Let $G$ be a connected graph. Then there is no subgraph of $G$ with the $Q$-spectrum  identical to ${\rm Spec}_{Q}(G) \cup \left\{ {{{\left[ 2 \right]}^1}} \right\}$. Moreover, If $G$ is of order at least 3, then $q_1(G)\geq 3$.
\end{lemma}
\begin{proof}
Suppose by the contrary that there is a subgraph of $G$, say $G'$, such that ${\rm Spec}_{Q}(G')={\rm Spec}_{Q}(G) \cup \left\{ {{{\left[ 2 \right]}^1}} \right\}$. But, in this case $ |E(G')|=|E(G)|+1$ and $ |V(G')|=|V(G)|+1$. Therefore there exists a vertex $v$ of $G'$ with the degree one such that
$G'-v=G$. This means that  $G$  is a proper subgraph of the connected graph $G'$ and so by Lemma \ref{lem:subg-q1}, $ q_1(G')> q_1(G)$, a contradiction. If $G$ is a connected graph of order at least 3, it has $K_3$ or $K_{1,2}$  as  its subgraph. Moreover, ${\rm Spec}_{Q}(K_3)=\{[4],[1]^2\}$ and
${\rm Spec}_{Q}(K_{1,2})=\{[3],[1],[0]\}$. Therefore by Lemma \ref{lem:q1GH}, $q_1(G)\geq 3$.
\end{proof}

\begin{theorem}\label{thm:GConDQS}
Let $G$ be a connected non-bipartite graph with $n\geq 3$ vertices which is {\rm DQS}. Then for any positive integer $r$, $G \sqcup rK_2$  is {\rm DQS}.
\end{theorem}

\begin{proof}
Let $H$ be a graph  $Q$-cospectral with $G \sqcup rK_2$. Then
by Lemmas \ref {lem 2-1} and \ref{lem:gipower}, $H$ has $n + 2r$ vertices, $n + 1+r$ edges and exactly $r$ bipartite components.
We perform mathematical induction on $r$. Suppose that $H$ is a graph  $Q$-cospectral with $G \sqcup K_2$. Then
 $${\rm Spec}_{Q}(H)={\rm Spec}_{Q}(G) \cup {\rm Spec}_{Q}(K_2)={\rm Spec}_{Q}(G) \cup\left\{ {{{\left[ 2 \right]}^1}}, [0]^1 \right\}.$$
Since  $G$  is a connected non-bipartite graph, by Lemma \ref{lem 2-1}, it has not $0$ as its signless Laplacian eigenvalue.
Therefore, $H$  has exactly one bipartite component. Therefore, by Lemma \ref{lem:thm4} we get $H=G\sqcup K_2$.
Now, let the assertion holds for $r$; that is, if ${\rm Spec}_{Q}(G_1)={\rm Spec}_{Q}(G) \cup {\rm Spec}_{Q}(rK_2)$, then $ G_1=G\sqcup rK_2 $. We show that it follows from ${\rm Spec}_{Q}(K)={\rm Spec}_{Q}(G) \cup {\rm Spec}_{Q}((r+1)K_2)$ that $ K=G\sqcup (r+1)K_2$. Obviously, $ K $ has $2$ vertices, one edge and one bipartite component more than $ G_1 $. So, we must have $ K=G_1\sqcup K_2 $. Now, the inductive hypothesis holds the proof.
\end{proof}

Lemma \ref{lem:RegularCone} and Theorem \ref{thm:GConDQS} imply the following corollary.

\begin{corollary}
  For a $k$-regular graph $G$ of order $n$,
$(G\nabla K_1)\sqcup rK_2$ is   {\rm DQS}   if either of the following conditions holds:
\begin{enumerate}[{\rm(i)}]
  \item $k\in\{1, n-2\}$,
  \item $k = 2$ and $n\ge 11$,
  \item $k = n-3$ and $\overline{G}$ has no triangles.
\end{enumerate}
\end{corollary}

Lemma \ref{lem:RegularConeK2} and Theorem \ref{thm:GConDQS} imply the following corollary.

\begin{corollary}\label{cor:RegularConeK2}
Let $G$ be a $k$-regular graph of order $n$. Then $(G\nabla K_2)\sqcup rK_2$ is
 {\rm DQS} for $k\in\{1,n-2\}$. For $k= n-3$,
$(G\nabla K_2)\sqcup rK_2$ is  {\rm DQS} if $\overline{G}$ has no triangles.
\end{corollary}


Lemma \ref{lem:KnMatching} and Theorem \ref{thm:GConDQS} imply the following corollary.

\begin{corollary}\label{cor 3-2}
 Let  $G$  be a non-bipartite graph obtained from $K_n$ by deleting a matching. Then
 $ G\sqcup rK_2 $  is {\rm DQS}.
 \end{corollary}

\begin{remark} Some $3$-regular  {\rm DAS}  graphs are given in \cite{HuangK} and references therein. We can obtain {\rm DQS}   graphs with independent edges from Corollary \ref{cor 3-2}.
\end{remark}

Lemmas \ref{lem:regularDAQS} and \ref{lem:regularDAS} and Theorem \ref{thm:GConDQS} imply the following corollary.
\begin{corollary}
Let $G$ be a  $k$-regular connected non-bipartite graph of order $n$. Then $ G\sqcup rK_2 $ is {\rm DQS} if
either of the following holds
\begin{enumerate}[{\rm(i)}]
\item
$k\in \{2, n-1, n-2, n-3\}$.
\item $k=n-4$ and $G$ is  {\rm DAS}.
\end{enumerate}
\end{corollary}


Lemma \ref{lem:main}  and Theorem  \ref{thm:GConDQS} imply the following corollary.

\begin{corollary}
Let $G$ be any of the following graphs. Then $G\sqcup rK_2$ is   {\rm DQS}.
\begin{enumerate}[{\rm(i)}]
\item
 The graphs $C_n$ ($n$ is odd), $K_n$ ($n\ge 4$).
 \item
 The graphs $\overline{P_n}$ ($n\ge 5$).
 \item The wheel graph $K_1\nabla C_n$.
\item Every lollipop graph $H_{n,p}$ when $p$ is odd and $n\ge 8$.
\item
The kite graph $Ki_{n,n-1}$ for $n\ge 4$ and $n\neq 5$.
\item The friendship graph $ F_n $.
\item
$(C_n\circ tK_1)$, when $n$ is odd and $(n \neq 32, 64)$  and $t\in\{1,2\}$.

\item $U_{r,n-r}$  if $r(\ge 3)$ is odd and $n\ge 7$.
\item
$CS(n,\alpha)$ when
$1\leq\alpha\leq n-1$ and $\alpha\neq 3$.
\item
$S(n,c)$ and its complement where $n\ge 2c+1$ and $c\ge 1$.
\item
$H\nabla K_m$ where $H$ is an
$(n-2)$-regular graph on $n$ vertices, and $\overline{K_n}\nabla K_2$ except for $n=3$.
\item
The dumbbell graphs $D_{p,k,q}$ ($p$ or $q$ is odd) different from $D_{3q,0,q}$ and all non-bipartite theta graphs $\Theta_{r,s,t}$.

\end{enumerate}
\end{corollary}


\end{document}